\documentclass[12pt]{paper}
\usepackage{amssymb,amsfonts,amsthm,amsmath}
\usepackage{latexsym,graphicx}

\usepackage{graphicx, verbatim}
\usepackage{xcolor}
\usepackage{verbatim}

\usepackage[a4paper, hmargin=3cm, vmargin={4cm, 4cm}]{geometry}

\usepackage{fancyhdr}
\usepackage{amssymb,amsfonts,amsthm,amsmath}
\usepackage{latexsym,multicol,graphicx}
\usepackage{float}
\usepackage{pgfplots}

\newtheorem*{LT}{Lindel\"of Theorem}
\newtheorem{corollary*}{Corollary}

\newtheorem*{def*}{Definition}

\newtheorem{theorem}{Theorem}
[section]

\newcommand{\eps}{\varepsilon}

\usepackage{ulem}
\theoremstyle{plain}
\pgfplotsset{compat=1.17} 
\pgfplotsset{soldot/.style={only marks,mark=*, line width=0.2pt, mark size=1pt}}
\newtheorem*{Lunc}{Lunc Theorem}
\newtheorem*{GS}{Grishin-Sodin Theorem}

\title{Uniqueness sets with angular density for spaces of entire functions, II: $\rho=1$. }
\author{ Anna Kononova
\thanks{Supported by  Israel Science Foundation grant  No. 1288/21, and by the Center for Integration in Science of the Israel's Ministry of Aliyah and Integration.}
}
\date{\today}
\begin{document}
\date{\today}
\maketitle

\begin{abstract}
    In this note, which is {the second part} of a 
three-part series,   we focus on uniqueness sets specifically in the case of spaces of entire functions of exponential type. As in the first part, we consider sets with angular density; however, now we abandon the Lindel\"of condition restriction.

{This part is essentially self-contained and can be read independently of the first one.}
\end{abstract}

\section{Introduction and the main result}

Let $\mathcal E_{1,\sigma
}$ be the set of the entire functions of exponential  type $\sigma$, that is, {for each $f\in \mathcal E_{1,\sigma}$ and for each $\eps>0$} there exists $C_\eps$ such that 
$$|f(z)|\le C_\eps e^{(\sigma+\eps)|z|}, \;\;\; z\in\mathbb C.$$

    For a given sequence $\Lambda=\{\lambda_k\in\mathbb C: 1\le|\lambda_1|\le|\lambda_2|\le\cdots<\infty\}$ without finite limit points, we will use the following notation for the number of points of $\Lambda$ in a sector:
    $$n_\Lambda(R;\alpha, \beta)=\#\{\lambda_k\in\Lambda:|\lambda_k|<R, \arg \lambda_k\in(\alpha,\beta)\}.$$

We say that the set $\Lambda$ has an {
\bf angular density} (with respect to the order $\rho=1$), if for all $0\le\alpha<\beta<2\pi,$ except, maybe, countably many values, there exists the limit
\begin{equation}
\label{DeltaL}
\Delta_\Lambda(\alpha, \beta)=\lim_{R\to\infty}\frac{n_\Lambda(R;\alpha,\beta)}{R},
\end{equation}
and we will consider $\Delta_\Lambda$ as a non-negative measure on $[0,2\pi)$.
{ The corresponding countable set of the exceptional values of arguments will be denoted by $X$, so that the condition \eqref{DeltaL} is satisfied for all $\alpha, \beta\notin X.$}

A set $\Lambda$ is called {\bf regular} if 
\begin{description}
    \item [ (i)] it has an angular density;
    \item[(ii)] the {\bf Lindel\"of condition} holds:
    there exists a finite limit 
    \begin{equation*}
    \lim_{R\to\infty}\sum_{0<|\lambda|\le R}\frac1{\lambda}.\end{equation*}
\end{description}

The regularity of the set $\Lambda$ implies (see \cite[Chapter II, Sec.3]{Levin}) that the measure $\Delta_\Lambda$ has zero first moment:
\begin{equation}\label{Lindelof}
    \int_0^{2\pi} e^{it}d\Delta_\Lambda(t)= 0.
\end{equation}

We will say that a   sequence  $\Lambda$ is  a {\bf uniqueness set} for $\mathcal E_{1,\sigma}$, if the only function in $\mathcal E_{1,\sigma}$ that vanishes at $\Lambda$ is the zero function:
$$f\in\mathcal E_{1,\sigma}\;\;\;{\it and}\;\;f|_\Lambda \equiv 0\;\;\Rightarrow \;\;\;f\equiv 0.$$ Otherwise, we will say that it is a     {\bf nonuniqueness set} for $\mathcal E_{1,\sigma}$.

{\it
Given  a discrete set $\Lambda\subset \mathbb C$, we let
\begin{align*}
    \sigma_U(\Lambda)
    =&\sup\{\sigma:\Lambda \text{ \;is\; a \;uniqueness\; set \;for\; }\mathcal E_{1, \sigma}\}\\
    =&\inf\{\sigma: \exists f\in \mathcal E_{1, \sigma}\setminus\{0\}: \;f|_\Lambda=0
    \}.
    \end{align*}
We call this quantity the {\bf critical uniqueness type of $\Lambda$}.

Analogously, we define the {\bf critical zero set type} $\sigma_Z(\Lambda)$ by
$$\sigma_Z(\Lambda)=\inf\left\{\sigma: \exists f\in \mathcal E_{1,\sigma}:\;\;f^{-1}\{0\}=\Lambda\right\}.$$
}


{ A $2\pi-$periodic function $h:\mathbb R\to\mathbb R$ is called  {\bf trigonometrically convex } if  for any  $t_1<t_2<t_3,$ such that  $t_3-t_1<\pi$ the following inequality holds:
$$ h(t_1)\sin (t_2-t_3)+h(t_2)\sin (t_3-t_1)+h(t_3)\sin (t_1-t_2)\le 0.
$$
}
Let us review some   key facts we will use (for other properties of trigonometrically convex functions see \cite[Chapter I, Sec.16]{Levin} and \cite[Lecture 8]{LevinLectures}).
\paragraph{Auxiliary facts.}
\paragraph{1.} If $f\in\mathcal E_{1,\sigma},$ then its indicator   $$h_f(t):= \limsup_{r\to\infty}\frac{\log|f(re^{it})|}{r}$$
is trigonometrically convex function. Moreover, $\displaystyle\max_{t\in[0;2\pi]}h_f(t)= \sigma.$
 \paragraph{2.} {Given a  trigonometrically convex function $h$ 
  the corresponding convex compact set $I_h\subset\mathbb C$ is defined by
 $$I_h:=\bigcap_{t\in[0,2\pi]} \{(x,y):x\cos t+y\sin t\le h(t)\}.$$
 The set 
 \begin{equation}
     \label{SuppLine}
 L_t=\{(x,y):x\cos t+y\sin t= h(t)\}
 \end{equation}
 is called  a {\bf support line} of the convex set $I_h.$ 
 } The correspondence between  $I_h$ and $h$ is  one-one. {Moreover, if $h_2(t)-h_1(t)=a\cos t+ b\sin t, a,b\in\mathbb R,$
 then $I_{h_2}$ can be obtained from $I_{h_1}$ by a translation by the vector $(a,b)$.} Note that in terms of indicators of entire functions the addition of the function $a\cos t+ b\sin t$ to the indicator $h_f$ corresponds to the multiplication of $f$ by factor $e^{(a-ib)z}$, so that the zero set of the function $f(z)e^{(a-ib)z}$ is the same as the zero set of $f$.

 \paragraph{3.} {Given a trigonometrically convex function $h,$   the  measure $\Delta_h$ is defined as follows:}
\begin{equation}
    \label{Delta}
2\pi\Delta_h(\alpha,\beta)=h'_+(\beta)-h'_-(\alpha)+\int_\alpha^\beta h.\end{equation}
 {Note
that if $h_2(t)-h_1(t)=a\cos t +b\sin t$, then 
$\Delta_{h_1}=\Delta_{h_2}$.}

 \paragraph{4.}
 From the geometrical point of view the value $\Delta_h([\alpha,\beta])$ is equal to the length of the corresponding arc of the set $\partial I_h$ divided by $2\pi.$ { This arc is defined as the intersection of the set $\partial I_h$ with all supporting lines $L_t$ for $t\in[\alpha,\beta]$.}

 \paragraph{5.} 
  Given a  measure $\Delta$ satisfying the moment condition \eqref{Lindelof}, the function 
 {\begin{equation}
    \label{h}\widetilde h_\Delta(t):=2\pi\int_{0}^{t}\sin (t-\varphi){\rm d}\Delta(\varphi)\end{equation}} is trigonometrically convex. Moreover, for $h=\widetilde h_\Delta(t)$ and $\Delta_h=\Delta$, relation \eqref{Delta} is satisfied \cite[Chapter I, Sec.19]{Levin}.


Given a bounded subset of the complex plane $E\subset\mathbb C$ we will denote by $D(E)$ its {\bf circumcircle}, which is defined to be the smallest circle containing the set $E$. The radius of this circle $R(E)$ is called the {\bf circumradius} and its center $O(E)$ is called the {\bf circumcenter}.

Given a regular set $\Lambda$ with angular density $\Delta$, put
 \begin{equation*}   
h_{\Delta}(t):=\widetilde h_{\Delta}(t) +A\cos t+B\sin t,\end{equation*}
where $\widetilde h_\Delta$ is defined by \eqref{h} and the parameters $A,B\in\mathbb R$ are chosen in such a way that the circumcenter $O(I_{h_\Delta})$ is located at the origin.
We will use the following notation: $D_\Delta:=D(I_{h_\Delta})$, $R_\Delta:=R(I_{h_\Delta})$.

{It is well known that, by the Levin-Pfluger theory (see \cite{Levin, AD1}),  the critical uniqueness type   of a  regular set $\Lambda$
  equals  the critical zero type $\sigma_Z(\Lambda)$.} Furthermore (see  \cite[Theorem 0.1]{Azarin-Giner}, \cite[Theorem 1.1]{Khabibullin2}), $\sigma_Z(\Lambda)$ is equal to the circumradius of the corresponding convex set:
$$\sigma_U(\Lambda)=\sigma_Z(\Lambda)=R_{\Delta_\Lambda}.$$

The same question for the non-regular sets, as far as the author knows, remains open.
The peculiarity of regular sets is connected with the classical Lindel\"of theorem (we provide here the formulation only for the case of  functions of the exponential type).

\begin{LT}{The entire function $f$  of exponential type with zero set $\mathcal Z$   belongs to the set $\mathcal E_{1,\sigma}$ if and only if }
\begin{enumerate}
    \item $n_{\mathcal Z}(R)\le CR,$ for some $C>0;$ 
    \item  the following sums are bounded:
    $$s(R)=\sum_{0<|z|\le R,  z\in\mathcal Z}\frac 1{z}.$$
\end{enumerate}

\end{LT}

{In particular, this theorem blocks  a non-regular sequence $\Lambda$  
  with angular density $\Delta_\Lambda$, such that the second condition of the Lindel\"of theorem is not satisfied,
  from being an exact zero set of any function of exponential type.}
On the other hand, one can always add an extra set of points to get a new set 
$\widetilde \Lambda\supset\Lambda$, so that there exists an
entire function of exponential  type, with zeros exactly at the points of the set   $\widetilde \Lambda$ (for example, one can consider the set $\widetilde \Lambda:=\Lambda\cup\{-\lambda_k,\lambda_k\in\Lambda$\} --- obviously,  this set satisfies the conditions of the Lindel\"of theorem). 

The main goal of this note is to find the value of the critical type $\sigma_U(\Lambda)$ in case when the  set $
\Lambda$ with angular density $\Delta_\Lambda$ is non-regular.
As we will show, for any non-regular set $\Lambda$ with angular density $\Delta_\Lambda$ there exists its "minimal regularization" $\Lambda^*$ such that
\begin{enumerate}
    \item $\Lambda^*\supset\Lambda$;
    \item $\Lambda^*$ is regular;
    \item $\sigma_U(\Lambda)=\sigma_Z(\Lambda^*)=\sigma_U(\Lambda^*).$
\end{enumerate}
Moreover, the minimal regularization can always be obtained from $\Lambda$ by adding two  sequences, one of which lies on a single ray, and the other has zero angular density.

The main result of this note is the following theorem.
\begin{theorem}\label{main}
For any discrete set $\Lambda$ with angular density $\Delta_\Lambda=\Delta$ we have
$$\sigma_U(\Lambda)=R_{ \Delta^*},$$
where
    $$\Delta^*:=\Delta +A_0\delta_{\alpha_0},\;\;\;A_0e^{i\alpha_0}=-\int_0^{2\pi} e^{it}d\Delta(t), \;\;A_0\ge 0.$$

\end{theorem}

\paragraph{Acknowledgements}
The author sincerely thanks Alexander Borichev and Mikhail Sodin for their attention to this work, and for their valuable remarks.

The author also expresses gratitude for the support provided  by Israel Science Foundation grant  No. 1288/21, and by the Center for Integration in Science of the Israel's Ministry of Aliyah and Integration.

\section{Proof of the Theorem \ref{main}}

\subsection {$\sigma_U(\Lambda)\le R_{ \Delta^*}$}

Put 
$$\widetilde\Lambda:=\begin{cases}
    \displaystyle\Lambda\cup\left\{\frac {ke^{i\alpha_0}}{A_0},\;\;k\in\mathbb N\right\},&A_0>0,\\
    \Lambda,&A_0=0.
\end{cases}
$$
Clearly, $\Delta_{\widetilde\Lambda}=\Delta^*$ and since for $\Delta^*$ the moment condition \eqref{Lindelof} is fullfilled, the following  Lunc theorem can be applied. 

\begin{Lunc}{\rm \cite{Krivosh,Grishin,Lunc}}
    For any  discrete set $\Gamma$ with angular distribution $\Delta_{\Gamma}$ with zero first moment
    there exists a set $\Gamma^*\supset\Gamma,$ such that 
$\Delta_{\Gamma^*}=\Delta_{\Gamma}$   and 
\begin{equation*}
\lim_{r\to\infty}\left(\sum_{|\lambda_k|<r, \lambda_k\in\Gamma^*}\frac1{\lambda_k}\right)=0.
    \end{equation*}
\end{Lunc}

{The proof of this theorem can be found in various sources
 (see \cite{Krivosh}, \cite {Grishin}, \cite{Lunc}). 
 For instance, in the most recent reference \cite{Krivosh}, where a more general result is established, the Lunc theorem is a special case of  \cite[Theorem 2.9]{Krivosh}. 
 }

It follows from the  Lunc theorem that there exists  a set $\Lambda^0$ of zero angular density such that the set
$ \Lambda^*:=\widetilde\Lambda\cup\Lambda^0$ is regular. 

{Recall, that by the Levin-Pfluger theory \cite[Chapter II]{Levin} (see also \cite[Corollary 1.1]{AD1}) $\sigma_Z(\Lambda^*)=R_{\Delta^*}$, hence for any 
$\sigma\ge R_{\Delta^*}$   there exists  
$f\in\mathcal E_{1,\sigma}$ such that 
$f|_{\Lambda^*}=0$, and therefore $f|_{\Lambda}=0$.   As a result, we have 
$\sigma_U(\Lambda)\le R_{\Delta^*}.$
}

\subsection{$\sigma_U(\Lambda)\ge R_{\Delta^*}$}

In this section we are going to prove that any set 
$\Lambda$ with angular density $\Delta$ is a uniqueness set of the space $\mathcal E_{1,\sigma}, $ provided that  
$\sigma<R_{\Delta^*}.$

Suppose that there exists  $g\in\mathcal E_{1,\sigma}$ where $\sigma<R_{\Delta^*},$ and such that $\left.g\right|_\Lambda=0.$
We will use the following uniqueness theorem proved by Grishin and Sodin \cite{GS} to show that in such a case $g\equiv 0$:

\begin{GS} {\rm \cite[Theorem 6.2]{GS}}
Given a discrete set $\Lambda:=\{r_ke^{i\theta_k}\}$, let $F$ be an entire function of order $1$ such that $F|_\Lambda=0$.
If there exists a non-negative trigonometrically convex function $k$ such that
\begin{equation}
    \label{GSe}
\displaystyle\limsup_{R\to\infty}\frac 1{\log R}\int_1^R\frac{\sum_{r_n\le R}k(\theta_n)}{r^{2}}{\rm d}r>\int_0^{2\pi}\log|F(\varphi)|\;{\rm d}\Delta_k(\varphi),
\end{equation}
then $F\equiv 0.$
\end{GS}

{ To apply Grishin-Sodin Theorem we need to choose an appropriate non-negative trigonometrically convex function $k$.  An important consideration here is that, in order to use Grishin-Sodin Theorem, we need to have some information on the density of zeros (or at least on a lower bound for the density) of entire function $g$ with arguments $\theta$ where the function $k$
 is strictly positive. 
Note that in case $A_0\ne 0$ we know a priori that there is a "gap" in our information regarding the distribution of the zeros of the function $g$: all that we have is the distribution $\Delta$ of the non-regular part $\Lambda$ of the set of zeros $\mathcal Z_g$ of the function $g$ of  exponential type. Therefore, we strongly suspect that the angular distribution of the set $\mathcal Z_g$   in the case $A_0>0$ should contain something additional near the point
$\alpha_0$, but we cannot be sure that adding just one mass point to the measure will always be sufficient. Due to this lack of information, it seems preferable to choose a trigonometrically convex function $k$ that vanishes at the point $\alpha_0.$ }

The proof is divided into two parts: we start with the case of a finitely supported measure $\Delta$, then we
extend the result to the case of general measure.

\subsubsection{A finitely supported angular density.}

Assume that $$\Delta=\sum_{k=1}^N A_k\delta_{\alpha_k},$$ {where $0\le\alpha_0\le \alpha_1<\alpha_2<\dots<\alpha_N<2\pi.$
In the case $A_0=0$  we fix $\alpha_0=\alpha_1.$ }

Consider a new measure $\;\Delta^*=\sum_{k=0}^NA_k\delta_{\alpha_k}.$
Then $I_{\Delta^*}$ is a convex polygon with vertices $C_k\in\mathbb C,$ $ k=0,\dots, N,$ such that
$$C_{k+1}-C_k=2\pi i\cdot A_k\cdot e^{i\alpha_k}, $$
where $C_0=C_{N+1}$ (see \cite[Chapter I, \S 19]{Levin}). That is, the edge  $[C_kC_{k+1}]$  of the polygon  $I_{\Delta^*}$  is perpendicular to the ray with angle $\alpha_k$, and has the length $|C_kC_{k+1}|=2\pi A_k.$

Without loss of generality, we can assume that  $ D_{\Delta^*}$ (the circumcircle of $I_{\Delta^*}$) is centered at the origin. Put $M_{\Delta^*}=I_{\Delta^*}\cap \partial D_{\Delta^*}.$


 From the definition of the circumcircle it follows that there are two possibilities:
either there is a pair of diametrically opposite points $M_1, M_2\in M_{\Delta^*},$  or there are three points $M_1, M_2, M_3\in M_{\Delta^*}$  such that the origin is the inner point of the triangle $\mathcal T_M$ with vertices $ M_1, M_2, M_3.$

{ \paragraph{ Case 1.} Suppose that there are two diametrically opposite points $M_1=C_{j_1}, M_2=C_{j_2},\; j_1<j_2,$ such that $M_{1}=-M_{2}=e^{i\mu_{1}}R_{\Delta^*},
\;\;\mu_{2}=\mu_{1}+\pi.
$ Without loss of generality we can assume that $0\le\alpha_0\le\mu_1<\mu_2<2\pi.$

We put
$$k^*(t)=
\begin{cases}
    \sin(t-\mu_{1}),& t\in [\mu_{1},\mu_{2}],\\
    0, &{\rm otherwise}.
\end{cases}$$ 

Note that this function meets our preferences for the choice of a  trigonometrically convex function:  the function $k^*$ is nonnegative  and  $k^*(\alpha_0)=0$. 

Let $\eps>0$ (to be defined later) { be such that $\eps <A_k,\;k=1,\dots, N.$ }
Using the uniform continuity of the function $k^*$  we can choose $\delta>0$ so that for any $t_1, t_2 $ with $|t_1-t_2|<\delta$, the inequality $|k^*(t_1)-k^*(t_2)|<\eps$  holds. Note that, using the countability of the exceptional set $X$ in \eqref{DeltaL}, we can choose $\delta$ in such a way that $\forall j \;\;\alpha_j\pm \delta\notin X$.
Now, since for all $\alpha, \beta\in[0,2\pi)\setminus X$ we have 
$$\lim_{r\to\infty}\frac{n_\Lambda(r;\alpha, \beta)}{ r}=\Delta(\alpha,\beta),$$ 
there exists $R>0$ (depending on $\varepsilon$) such that  
$\forall j=1,\dots,N$ 
$$\left|\frac{n_\Lambda(r, \alpha_j-\delta,\alpha_j+\delta)}{r}-A_j\right|<\eps, \;\;\forall r>R.$$

We proceed by estimating the left-hand side of the inequality (\ref{GSe}) ({put $r_k:=|\lambda_k|$, $\theta_k:=\arg(\lambda_k)$}).
\begin{multline*}
\displaystyle\limsup_{R\to\infty}\frac 1{\log R}\int_1^R\frac{\sum_{r_n\le r}k^*(\theta_n)}{r^{2}}{\rm d}r\ge \sum_{k=1}^N (A_k-\eps)\; (k^*(\alpha_k)-\eps)\\\ge\sum_{\alpha_k\in[\mu_{1},\mu_{2}]} A_k\; \sin(\alpha_k-\mu_{1})-\eps\cdot\left(N+\sum_{k=1}^N A_k\right). 
    \end{multline*}

From the geometrical point of view the value of the expression $2\pi A_k\; \sin(\alpha_k-\mu_{1})$ is equal to the length of projection of the edge  $[C_kC_{k+1}]$ of the convex set $I_{\Delta^*}$ to the diameter $[M_{1}M_{2}], $ and the sum of all these projections over all ${\alpha_k\in[\mu_{1},\mu_{2}]}$  is equal to $ 2R_{\Delta^*}.$ We therefore obtain}
\begin{equation*}
\displaystyle\limsup_{R\to\infty}\frac 1{\log R}\int_1^R\frac{\sum_{r_n\le r}k^*(\theta_n)}{r^{2}}{\rm d}r\ge
 \frac1\pi R_{\Delta^*}-\varepsilon (N+ R_{\Delta^*}),
\end{equation*}
where in the last inequality we have used the fact that for the convex set $I$ the length of its boundary $\partial I$ is not greater then the length of the circumcircle $D(I),$
{that is,
$$2\pi \sum_{k=1}^N A_k\le 2\pi R_{\Delta^*}.$$
    }

Next, we  estimate the right-hand side of  inequality (\ref{GSe}).
    Since $$\Delta_{k^*}=\frac{1}{2\pi}(\delta_{\mu_1}+\delta_{\mu_2}),$$ we have
$$\int_0^{2\pi} h_g(\varphi){\rm d}\Delta_{k^*}(\varphi)\le \sigma\cdot \Delta_{k^*}([0,2\pi))=\frac\sigma\pi.$$
Under the assumption $\sigma<R_{\Delta^*},$ there exists $\eps>0$ such that
$$
\frac {R_{\Delta^*}}\pi- \varepsilon (N+ R_{\Delta^*})>\frac \sigma\pi,$$
hence, we obtain
$$
\displaystyle\limsup_{R\to\infty}\frac 1{\log R}\int_1^R\frac{\sum_{r_n\le R}k^*(\theta_n)}{r^{2}}{\rm d}r>\int_0^{2\pi} h_g(\varphi){\rm d}\Delta_{k^*}(\varphi),$$
and by Grishin-Sodin Theorem we get   $g\equiv 0.$

\paragraph{ Case 2.}

\begin{figure}[ht]
    \centering
\includegraphics[width=0.5\linewidth]{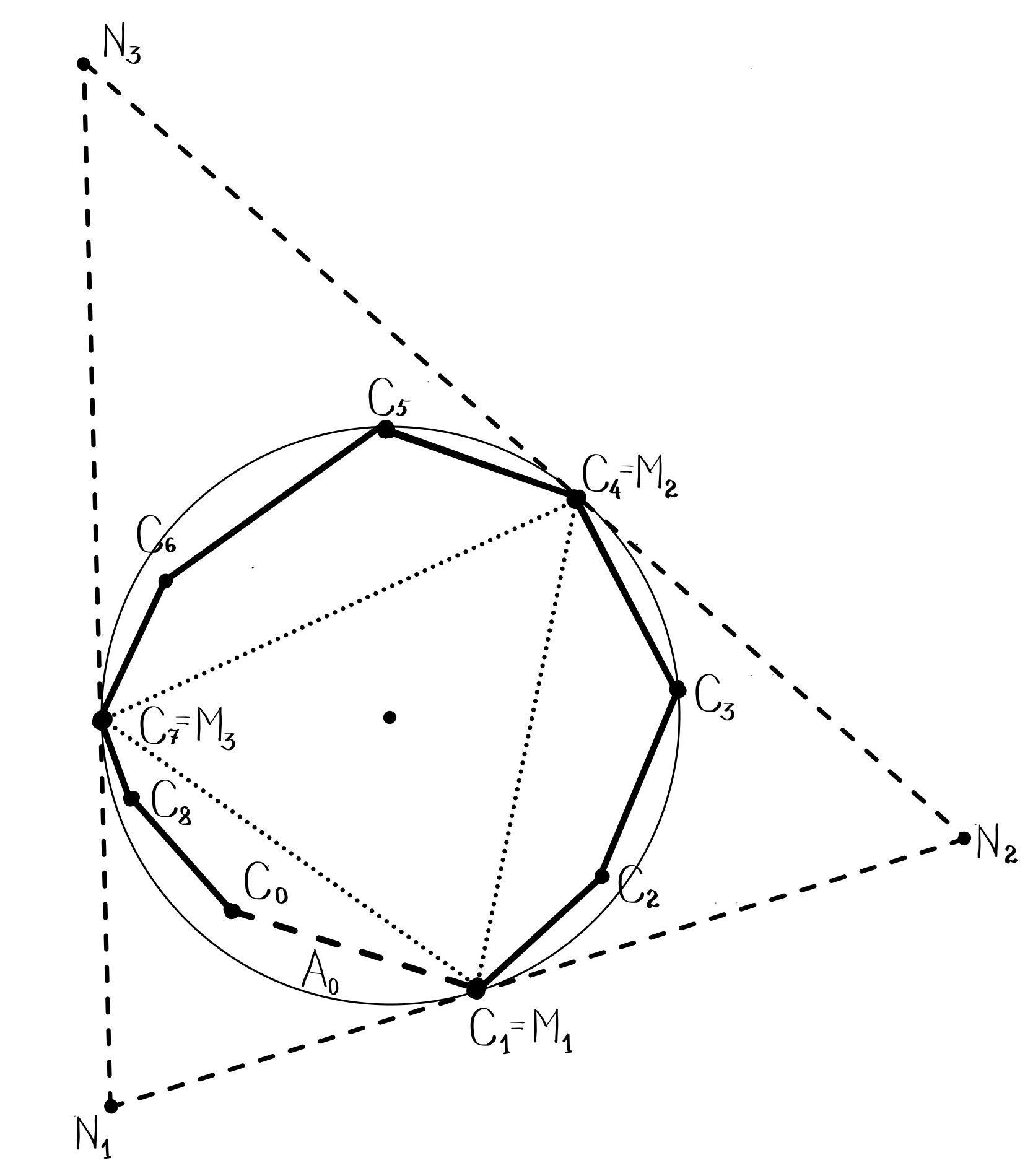}
    \caption{ Triangles $\mathcal T_M$ and $\mathcal T_N$.
}
    \label{fig1}
\end{figure}

Assume now that {there is no diameter 
 of $D_{\Delta^*}$ with endpoints in the set $M_{\Delta^*}$.}
It then follows   that there exists a triple of points $\{M_{1}, M_{2}, M_3\}\subset M_{\Delta^*}$: 
$$M_j=C_{k_j}=e^{i\mu_j}R_{\Delta^*}, \;\;j=1,2,3,$$
such that the point $0$ is an interior point of the triangle $\mathcal T_M$ with vertices $ M_1, M_2, M_3.$  
We assume without loss of generality, that $0\le\alpha_0\le\mu_1<\mu_2<\mu_3<2\pi.$
The polar set (with respect to $\partial D_{\Delta^*}$) to the triangle $\mathcal T_M$ is another triangle $\mathcal T_N$, which is tangential to the circle $D_{\Delta^*}$ (see Fig.\ref{fig1}). We denote its vertices as $ N_1,N_2,N_3,$ and suppose that the vertex $N_1$ is dual to the edge $[M_3M_1]$ (this edge is of especial importance for us, as it corresponds to the interval $[\mu_3,\mu_1]$ which contains the point~$\alpha_0$).

In this case we choose $k^*$ to be the support function of the triangle $\mathcal T_N$ viewed from the point $N_1$:
$$k^*(\theta)=\max_{j=1,2,3}\{ (|N_1N_j|\cdot \cos(\theta-\gamma_j))\}, $$
where,  $N_k-N_1=|N_1N_k|e^{i\gamma_j}.$ 
Note that, since for $k=1$ we have $|N_1N_k|=0$, this function is non-negative, as required by Grishin-Sodin Theorem. { Additionally, in accordance with our preferences for the choice of  $k^*$, it vanishes in the interval $[\mu_3,\mu_1]$, in particular $k^*(\alpha_0)=0$.} 
{Note also that,  due to the construction, we have $\gamma_2=\mu_1+\pi/2, \gamma_3=\mu_3-\pi/2,$ therefore, we have
\begin{gather}\label{k^*}
     k^*(\theta)=
\begin{cases}
    |N_1N_2|\cdot \sin(\theta-\mu_1),& \theta\in[\mu_1,\mu_2];\\
|N_1N_3|\cdot \sin(\mu_3-\theta)),&\theta\in[\mu_2,\mu_3];\\
0,&\theta\notin[\mu_1,\mu_3].
\end{cases}
\end{gather}
}

As in the previous case, we fix $\eps>0$ such that $\eps<A_k,\; k=1,\dots, N$, and choose $\delta>0$ such that for any $t_1, t_2 $ with $|t_1-t_2|<\delta$ it holds $|k^*(t_1)-k^*(t_2)|<\eps$, and such  that $\forall j \;\;\alpha_j\pm \delta\notin X$.
Then there exists $R>0$ (depending on $\varepsilon$) such that  
$\forall j=1,\cdots,N$ 
$$\left|\frac{n_\Lambda(r, \alpha_j-\delta,\alpha_j+\delta)}{r}-A_j\right|<\eps, \;\;\forall r>R.$$
Hence, analogously to the Case 1, we get
\begin{multline}\label{GS}
\displaystyle\limsup_{R\to\infty}\frac 1{\log R}\int_1^R\frac{\sum_{r_n\le r}k^*(\theta_n)}{r^{2}}{\rm d}r\ge \sum_{k=1}^N (A_k-\eps)\; (k^*(\alpha_k)-\eps)\ge\\\sum_{k=1}^N A_k\; k^*(\alpha_k)-\eps(N+R_{\Delta^*}).
    \end{multline}

{  Referring to the geometric interpretation, we now will  show that }
$$\pi\displaystyle\sum_{k=1}^N A_k\; k^*(\alpha_k)= S(\mathcal T_N),$$
where $S(\mathcal T_N)$ denotes the area of the triangle $\mathcal T_N.$ 

\begin{figure}[ht]
    \centering
\includegraphics[width=0.5\linewidth]{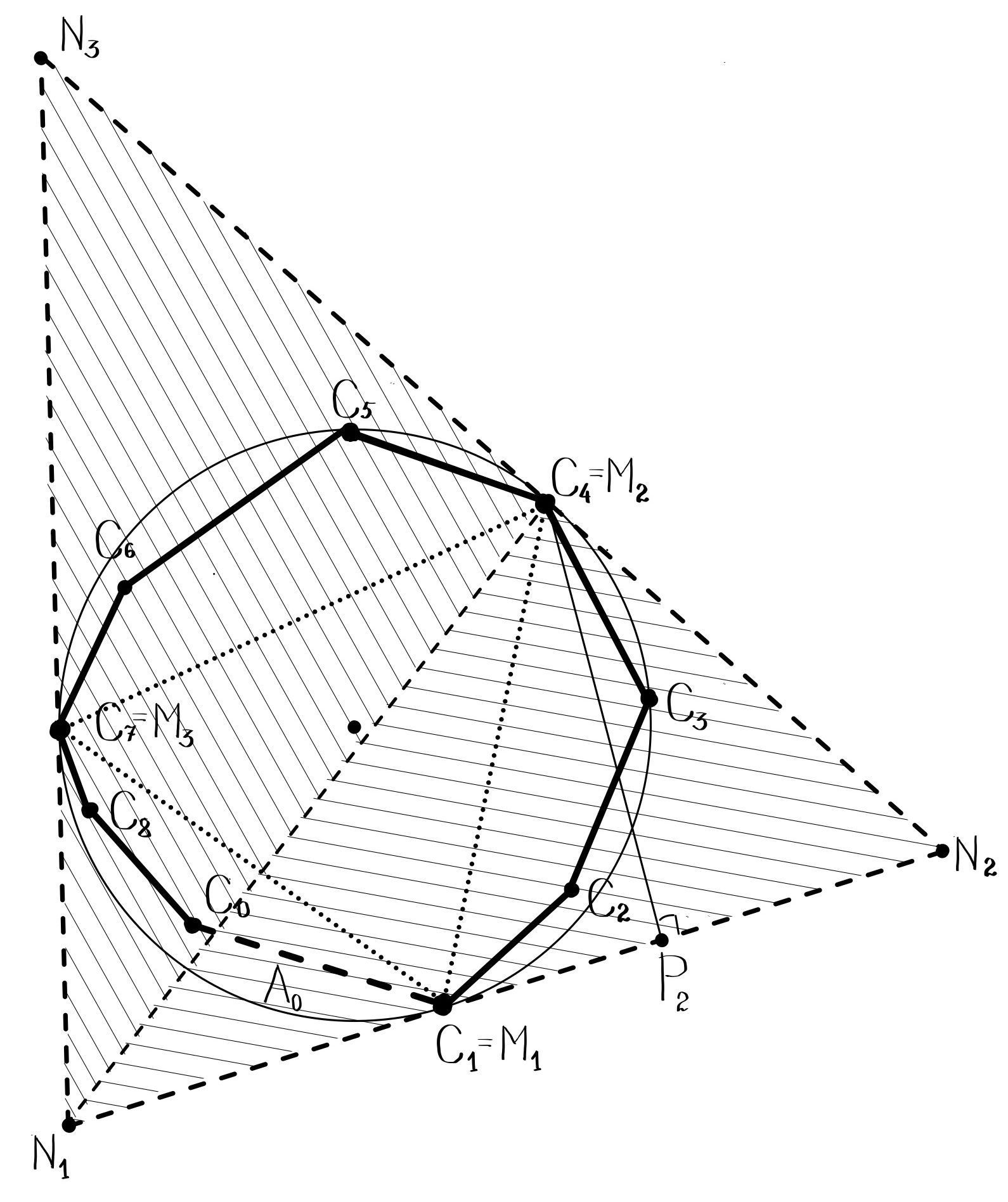}
    \caption{ \small 
    $S(\mathcal T_N)=S(N_1N_2M_2)+S(N_1N_3M_2)$, where
    $S(N_1N_2M_2)=\frac12|N_1N_2|\cdot|M_2P_2|$, 
   and similarly for the second triangle.}
    \label{fig2}
\end{figure}

{Recall that, by \eqref{k^*}, for $\theta=\alpha_k$ with $k\in[k_{1},k_2-1]$, we have 
$$2\pi A_k\cdot k^*(\alpha_k)=|C_kC_{k+1}|\cdot   |N_1N_2|\cdot \sin(\theta-\mu_1)=|N_1N_2|\cdot {\rm Proj}_2(C_kC_{k+1}), $$
where $ {\rm Proj}_2(C_kC_{k+1})$ is the length of the orthogonal projection of the edge $[C_kC_{k+1}]$ to the direction orthogonal to the line $(N_1N_2).$
Note that the length of the projection of the entire polygonal chain $C_{k_1}C_{k_1+1}\cdots C_{k_2}$, which connects the points $M_1=C_{k_1}$ and $M_2=C_{k_1+1}$, is equal to the length of the altitude $[M_2P_2]$ of the triangle $N_1N_2M_2$ from the vertex $M_2$ (see Fig. \ref{fig2}). Hence,
$$
   2\pi \sum_{k=k_{1}}^{k_{2}-1} A_kk^*(\alpha_k)=|N_1N_2|\cdot  |M_2P_2|.
$$
Analogously, for $k\in[k_{2},k_3-1]$, we have 
$$2\pi A_k\cdot k^*(\alpha_k)=|C_kC_{k+1}|\cdot   |N_1N_3|\cdot \sin(\mu_3-\theta)=|N_1N_3|\cdot {\rm Proj}_3(C_kC_{k+1}), $$
where $ {\rm Proj}_3(C_kC_{k+1})$ is the length of the orthogonal projection of the edge $[C_kC_{k+1}]$ to the direction orthogonal to the line $(N_1N_3).$ Taking the sum over all $k\in[k_{2},k_3-1]$, we get
$$
   2\pi \sum_{k=k_{2}}^{k_{3}-1} A_kk^*(\alpha_k)=|N_1N_3|\cdot  |M_2P_3|,
$$
where $ P_3$ is the foot of the altitude   of the triangle $N_1N_3M_2$ from the vertex $M_2$.
Hence, we obtain
\begin{equation}
\label{Ss}
   \pi \sum_{k=1}^{n} A_kk^*(\alpha_k)=\frac{1}{2}(|N_1N_2|\cdot  |M_2P_2|+|N_1N_3|\cdot  |M_2P_3|)=S(\mathcal T_N).
\end{equation}

}

On the other hand, 
$$\int_0^{2\pi} h_g(\varphi){\rm d}\Delta_{k^*}(\varphi)\le \sigma\cdot \Delta_{k^*}([0,2\pi)),$$
Recall that by the auxiliary fact 4 we have  (here, $N_4=N_1$)
$$\Delta_{k^*}([0,2\pi))=\frac1{2\pi}\cdot\sum_{j=1}^3 |N_jN_{j+1}|.$$
Hence,

\begin{equation}\label{S1}
    \int_0^{2\pi} h_g(\varphi){\rm d}\Delta_{k^*}(\varphi)\le\frac\sigma{2\pi}\cdot\sum_{j=1}^3 |N_jN_{j+1}|= \frac {\sigma \;S(\mathcal T_N)}{\pi R_{\Delta^*}}.
    \end{equation}

From \eqref{Ss} and the assumption $\sigma<R_{\Delta^*},$ it follows that there exists $\eps>0$ such that
\begin{equation}
\label{S4}
    \sum_{k=1}^N A_k k^*(\alpha_k)- \eps (N+ R_{\Delta^*})=\frac{S(\mathcal T_N)}{\pi}- \eps(N+ R_{\Delta^*})>\frac {\sigma\;S(\mathcal T_N)}{\pi R_{\Delta^*}}.
\end{equation}
Thus, combining \eqref{S4} and \eqref{S1}, we obtain
$$
\displaystyle\limsup_{R\to\infty}\frac 1{\log R}\int_1^R\frac{\sum_{r_n\le R}k^*(\theta_n)}{r^{2}}{\rm d}r>\int_0^{2\pi} h_g(\varphi){\rm d}\Delta_{k^*}(\varphi).$$
By Grishin-Sodin Theorem it follows that   $g\equiv 0.$

\subsubsection{General case.}

Suppose now that $\Delta$ is a finite Borel measure, and $\Delta^*:=\Delta+A_0\delta_{\alpha_0}.$ 
As before, we denote by $I_{\Delta^*}$ the corresponding convex set,  by $D_{\Delta^*}$ its circumcircle with center at the origin, and by $R_{\Delta^*}$ the circumradius of the set $I_{\Delta^*}$.

Put 
$M_{\Delta^*}:= I_{\Delta^*}\cap \partial D_{\Delta^*}$. 
As in the case of finitely supported measures, there are two scenarios depending on whether there are two diametrically opposite points in $M_{\Delta^*}$ or not.

{ 
\paragraph{ Case 1.} Suppose that there are two diametrically opposite points $M_1, M_2\in M_{\Delta^*},$ such that $M_{1}=-M_{2}=e^{i\mu_{1}}R_{\Delta^*},
\;\;\mu_{2}=\mu_{1}+\pi.
$ Without loss of generality we can assume that $0\le\alpha_0\le\mu_1<\mu_2<2\pi.$ 
In this case  we define
$$k^*(t)=
\begin{cases}
    \sin(t-\mu_{1}),& t\in [\mu_{1},\mu_{2}],\\
    0, &{\rm otherwise}.
\end{cases}$$ 

\paragraph{ Case 2.}
 Suppose now that there are no two diametrically opposite points in  the set $M_{\Delta^*}$. Then, from the definition of circumcircle it follows that there exist three points
$$M_j=e^{i\mu_j}R_{\Delta^*}\in M_{\Delta^*}, \;\;j=1,2,3,$$
such that
 the origin belongs to the interior of the triangle $\mathcal T_M$ with vertices $M_1, M_2, M_3.$
 We  assume, as before,  that $0\le\alpha_0\le\mu_1<\mu_2<\mu_3<2\pi,$
and denote by $\mathcal T_N$ the polar set (with respect to $D_{\Delta^*}$) to $\mathcal T_M$.

 As before, we denote by $N_1$ the vertex of the angle of $\mathcal T_N$  subtended by the arc $[M_3,M_1].$ The other two vertices we denote by $N_2, N_3$ in the counterclockwise order  (see Fig. \ref{fig3}).

\begin{figure}[ht]
    \centering
\includegraphics[width=0.65\linewidth]{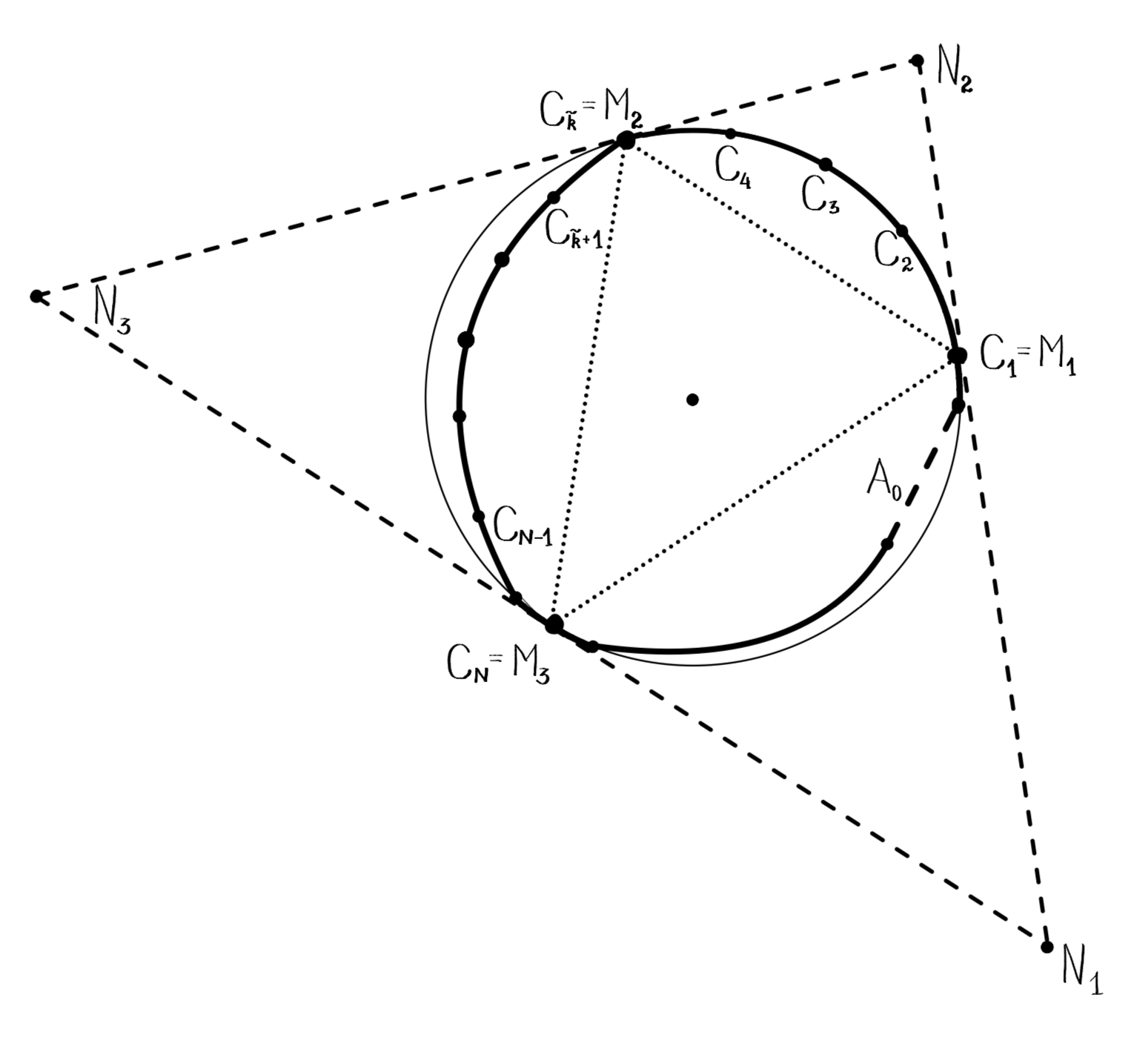}
    \caption{ Triangles $\mathcal T_M$ and $\mathcal T_N.$}
    \label{fig3}
\end{figure}

As in the case of discrete $\Delta$, in Case 2 we choose $k^*$ to be the support function of the convex set $\mathcal T_N$ viewed from the point $N_1$:
$$k^*(\theta)=\max_{j=1,2,3}\{{\rm Re} (|N_1N_j|\cdot e^{i(\gamma_j-\theta)})\}. $$
where, $N_j-N_1=|N_1N_j|e^{i\gamma_j}.$

\paragraph {Cases 1 $\&$ 2.}
In Case 1 we put $M_3:=M_2$ and $\mu_3:=\mu_2$.
Now, in both cases, in accordance with our preferences for the choice of  $k^*$, we obtained a function that vanishes on the interval $[\mu_3,\mu_1]$; in particular, $k^*(\alpha_0)=0$.
Now we will consider both cases Case 1 and Case 2 simultaneously.
 
{ Note that for every  point $M=R_{{\Delta^*}} e^{i\mu}\in M_{\Delta^*}$ the tangential line to the circumcircle $D_{\Delta^*}$  at this point serves as a support line to the set $I_{\Delta^*}$ at the point $M$. {More precisely,
$$L_{\mu}\cap D_{\Delta^*}=L_{\mu}\cap I_{\Delta^*}=\{M\},$$
where $L_{\mu}$ is the support line defined by \eqref{SuppLine}.}
Hence $\Delta^*(\{\mu\})=0.$ 
Therefore,  there are no positive masses of the measure $\Delta^*$ at the points   $\mu_1, \mu_2, \mu_3.$  Hence, {by \cite[Chapter III, Theorem 3]{Levin}}, $\mu_j \notin X$,$\;j=1,2,3$, where $X$ is the exceptional set.
 }

Fix $\eps\in(0,1) $ (to be defined later) and  choose $\delta>0$ such that for any $t_1, t_2$ with $|t_1-t_2|\le \delta$ we have $|k^*(t_1)-k^*(t_2)|<\eps.$

Let us choose $N\in\mathbb N$,  and a set of pairs $(\eta_k, C_k)$, where $\eta_k\in [0,2\pi)$ and  $C_k
\in \partial I_{\Delta^*}, k=1,\dots, N,$ in the following way.

 Put $C_1:=M_1, C_N:=M_3,$ and suppose that there exists $\widetilde k\in\{2,\dots, N\}$ such that $C_{\widetilde k}=M_2.$ 
To describe the requirements we impose on the rest of the points $C_k$,
note that each point $C_k$ belongs to the intersection of the set $\partial I_{\Delta^*}$ with some support line $L_{\eta_k}$, {the numbers $\eta_k$ to be defined later.} In the case when  $C_k$ is a corner point, there is an interval of possible values of $\eta_k$ (that is, the set of support lines contains more than one element). 

We impose the following conditions on our choice of $C_k$ and $\eta_k$
(we consider all the arguments as the numbers modulo $2\pi$):
\begin{enumerate}
  
    \item $ 0<\eta_{k+1}-\eta_{k}<\delta/3,\;\;\forall k=1,\cdots, N-1$;

    \item 
    $\eta_k\notin X\;\;k=1,\dots, N,$ where $X$ is the exceptional set (such a choice is always possible due to the countability of the set $X$).
\end{enumerate}

Now, we put 
$$C_0:=C_{N}$$
$$b_k:=C_{k+1}-C_k,\;\; \; B_k:=\frac{|b_k|}{ 2\pi}\ge0,\; k=0,\dots, N-1, $$
$$\beta_k=\begin{cases}
    \arg b_k-\pi/2,&b_k\ne 0;\\
    \dfrac{\eta_{k}+\eta_{k+1}}{2},&b_k= 0.
\end{cases}$$
{So, if $b_k\ne 0$, the value $\beta_k$ indicates a direction perpendicular to $[C_kC_{k+1}],$ while in case $b_k= 0$, the value $\beta_k$ corresponds to  the support line at the point $C_k=C_{k+1}.$
In any case, we have $\eta_k\le\beta_k\le\eta_{k+1}.$}
Analogously to the case of finitely supported angular density, we have
$$C_{k+1}-C_k=2\pi i\cdot B_k\cdot e^{i\beta_k}.$$
Thus, we obtained a convex polygon $C_1\dots C_N$. Note that for different values of $\eta_k$ our points $C_k$ may coincide. Hence, the case where some of  $B_j$ are equal to $0$ is not excluded.

Due to our construction, for each $t\in[\eta_k,\eta_{k+1}],$ $k=1,\dots, N-1,$ we have
$$|k^*(t)-k^*(\beta_k)|<\eps.$$

Now, 
there exists $R>0$ such that
$\forall k=1,\cdots,N-1$, for $r>R,$ we have 
$$\left|\frac{n_\Lambda(r, \eta_k,\eta_{k+1})}{r}-\Delta(\eta_k,\eta_{k+1})\right|<\eps/N.$$

 Then, for $r>R$, it holds
 {
\begin{multline*}
\displaystyle\sum_{r_n\le r}k^*(\theta_n)
\ge 
\sum_{j=1}^{N-1} n_\Lambda(r,\eta_j,\eta_{j+1}) (k^*(\beta_j)-\eps)\\
\ge
\sum_{j=1}^{N-1} \left(k^*(\beta_j)n_\Lambda(r,\eta_j,\eta_{j+1}) \; -\eps n_\Lambda(r,\eta_j,\eta_{j+1})\right){\rm d}r
\\\ge r\sum_{j=1}^{N-1} \left(k^*(\beta_j)\left(\Delta(\eta_j,\eta_{j+1})-\frac\eps N\right) -\eps\left(\Delta(\eta_j,\eta_{j+1})+\frac\eps N\right)\right).
\end{multline*}
}

Recall that
$2\pi \Delta(\eta_k,\eta_{k+1})$ is equal to the length of the corresponding arc of the boundary of the set $I_{\Delta^*}$. In our case it is the arc with the endpoints $C_k, C_{k+1},$ hence its length is not smaller than  $|C_kC_{k+1}|=2\pi B_k$. It follows  that 
$$
\Delta(\eta_k,\eta_{k+1})\ge B_k.
$$
Therefore, (put $K:=\displaystyle\max_{t\in[0,2\pi]} k^*(t)$) we get
{\begin{multline}\label{sumk*}
\displaystyle\frac1r{\sum_{r_n\le r}k^*(\theta_n)}
\ge 
\sum_{j=1}^{N-1} k^*(\beta_j)\Delta(\eta_k,\eta_{k+1})-K\eps\; -\eps\cdot\sum_{j=1}^{N-1} \Delta(\eta_k,\eta_{k+1})-\eps^2
\\\ge
\sum_{j=1}^{N-1} k^*(\beta_j)B_k\; -\eps\cdot(K+\Delta([0,2\pi])+1).
    \end{multline}
}

{

{Analogously to the case of a finitely supported  measure (we consider now the presence of degenerate edges does not affect the result), we have {$$\mathcal B:=\sum_{j=1}^{N-1}  B_j \cdot k^*(\beta_j)=\begin{cases}\dfrac{1}{\pi}R_{\Delta^*},& {\textrm{in case 1;}}
 \vspace{5 pt}\\
      \dfrac{1}{\pi}S(\mathcal T_N),& {\textrm {in case 2.}}
 \end{cases}$$}

Hence, using \eqref{sumk*}, we get
\begin{equation}
    \label{S}
\displaystyle\limsup_{R\to\infty}\frac 1{\log R}\int_1^R\frac{\sum_{r_n\le r}k^*(\theta_n)}{r^{2}}{\rm d}r\ge    \mathcal B - \eps\left( K+\Delta([0,2\pi))+1\right).  \end{equation}

On the other hand, as in the discrete case, we have
$$\int_0^{2\pi} h_g(\varphi){\rm d}\Delta_{k^*}(\varphi)\le \sigma\cdot \Delta_{k^*}([0,2\pi))= \frac {\sigma  \mathcal B}{ R_{\Delta^*}}.$$

Now,  for $\sigma<R_{\Delta^*}$ there exists $\eps\in(0,1)$ such that
$$
\mathcal B- \eps \left( K+\Delta([0,2\pi))+1\right)>\frac {\sigma \mathcal B}{ R_{\Delta^*}}.$$
Finally, from (\ref{S}) we get
$$
\displaystyle\limsup_{R\to\infty}\frac 1{\log R}\int_1^R\frac{\sum_{r_n\le R}k^*(\theta_n)}{r^{2}}{\rm d}r>\int_0^{2\pi} h_g(\varphi){\rm d}\Delta_{k^*}(\varphi).$$
 and by Grishin-Sodin Theorem it follows that  $g\equiv 0.$}
}

{\bf Remark.} Note that under the assumptions of the theorem, in the case $A_0>0,$ the additional measure $\Delta_0$, such that $\sigma_U(\Lambda)=R_{\Delta+\Delta_0}$, is not unique. This measure corresponds to the additional edge with the length $A_0.$ We can replace this edge by any curve such that the corresponding  set remains convex. 

\bigskip

\bigskip
\medskip

\noindent School of Mathematics, Tel Aviv University, Tel Aviv, Israel
\newline {\tt anya.kononova@gmail.com}

\end{document}